\documentclass{ article}

\usepackage{graphicx} 
\usepackage{amsmath,amssymb,dsfont,units,booktabs}
\usepackage{amsthm}
\newtheorem{theorem}{Theorem}[section]

\usepackage{natbib}
\usepackage{tabularx}
\usepackage{appendix}
\usepackage{graphicx, color}
\usepackage{color}
\usepackage[utf8]{inputenc} 
\usepackage{units}
\usepackage{siunitx}[=v2]
\usepackage{subcaption}
\usepackage{url}
\usepackage{soul}
\usepackage{algorithm}
\usepackage{amsthm}

\newcommand{\vt}{\mathbf{v}}
\newcommand{\xt}{\mathbf{x}}
\newcommand{\phit}{\mathbf{\phi}}

\newcommand{\Vt}{\mathbf{V}}

\renewcommand{\div}{\operatorname{div}}
\newcommand{\tr}{\operatorname{tr}}

\newcommand{\sigmat}{\boldsymbol{\sigma}}
\newcommand{\epsilont}{\dot{\boldsymbol{\epsilon}}}
\usepackage{bm}
\usepackage{algpseudocode}

\title{A Hybrid Particle-Continuum Method for Simulating Fast Ice via Subgrid Iceberg Interaction
}
\author{Carolin Mehlmann, Saskia Kahl}

\begin{document}

\maketitle

\begin{abstract}

A significant fraction (4\%-13\%) of Antarctic sea ice remains stationary as landfast sea-ice ("fast ice"), typically anchored by grounded icebergs. Current global climate models do not represent fast-ice formation due to iceberg grounding, as iceberg--sea-ice interaction mostly occurs at subgrid scales. We propose a novel subgrid-scale coupling mechanism between Lagrangian iceberg particles and an Eulerian sea-ice continuum model. This hybrid particle-continuum approach integrates feedback from icebergs into the sea-ice momentum equation via a Green’s function, a Stokeslet, representing the drag exerted by a point force on the viscous-plastic medium.
The coupled system, including the Stokeslet induced drag, is discretized using a finite-element method with piecewise linear basis functions.  The approach assumes that individual icebergs have diameters smaller than the grid spacing. The presented finite-element discretization is compatible with existing unstructured-mesh ocean model frameworks such as FESOM and ICON, ensuring practical applicability in Earth system modeling.
This work provides and analyzes, for the first time, a stable numerical framework to capture the effects of individual subgrid-scale icebergs on sea-ice dynamics. We derive an a-priori stability estimate bounding a functional of the sea-ice system and show that the momentum equation including the subgrid iceberg--sea-ice drag remains stable. Numerical test cases demonstrate the capability of the approach to capture fast-ice formation due to subgrid iceberg grounding on coarse horizontal grids.



\end{abstract}
\maketitle                   

\section{Introduction}
Around 4$\%$ to 13$\%$ of sea ice remains stationary, forming a narrow band around Antarctica, unlike the majority of sea ice, which drifts with winds and ocean currents ("pack ice"). This stationary landfast sea-ice ("fast ice") is anchored to the coastline or grounded by icebergs and has wide-ranging ramifications for global climate \citep{Fraser2023}. Additionally, {the fast-ice extent has strongly declined in recent years} \citep{Fraser2023}. The Intergovernmental Panel on Climate Change (IPCC) expresses low confidence in observed trends in fast-ice extent due to the limited observational record available  \citep{Meredith2019}. Further contributing to this uncertainty is the fact that current global climate models poorly represent Antarctic fast ice \citep{Fraser2023}, which casts doubt on their future projections for fast-ice and for sea-ice conditions around Antarctica in general \citep{Roach2020}.

In the Arctic, fast ice mainly forms through coastal attachment or grounding in shallow waters \citep{Fraser2023} and is modeled using basal stress, coastline drag, or rheology modifications \citep{KoenigBeatty_Holland2010,Lemieux2016,Liu2022}. Modeling fast ice in the Southern Ocean presents additional challenges, as mainly small grounded icebergs are the primary attachment mechanism in most areas \citep{Fraser2023, Fraser2012}.
Regional high-resolution sea-ice models ($\leq$ 2km resolution) can simulate fast ice by combining rheological modifications with static iceberg masks (e.g.,\citep{Huot2021, VanAchter2022}), but these approaches are computationally infeasible at global scales and fail to capture iceberg dynamics such as ungrounding. To address this scale mismatch, we propose a hybrid particle-continuum framework that enables the representation of dynamic grounding events by subgrid-scale icebergs in coarse-resolution models. The approach couples Lagrangian iceberg particles with an Eulerian continuum sea-ice model by introducing a subgrid-scale drag term. 

Existing sea-ice models that include iceberg--sea-ice interaction fall into two main categories: rheology-based and drag-based approaches. The former modifies the sea-ice constitutive relation to account for icebergs, typically on the scale of hundreds of meters, either resolved by the mesh or represented via local density functions \citep{Kahl2025,Vankova2017}. The drag-based approach \citep{Hunke2010} applies explicit interaction terms, but assumes large tabular icebergs ($\geq$40 km) resolved by the horizontal mesh, and deactivates drag when iceberg and sea-ice velocities align (due to numerical instabilities) or in grounded iceberg cases.
We extend this second class by introducing a subgrid-scale drag formulation based on Stokeslets, enabling icebergs smaller than the grid scale to exert localized forces on the sea-ice continuum. Our approach is stable and explicitly captures grounded iceberg forcing, making it well-suited for modeling Antarctic fast-ice dynamics at coarse resolution.

The core of this coupling mechanism is a Green’s function approach: each iceberg exerts a localized drag force on the sea-ice field, modeled using the Stokeslet \citep{Cortez2001}. This approach is related to force-coupling and immersed-boundary methods used in particle-fluid systems \citep{Lomholt2003,Peskin2002}, but is adapted here to model grounding events in a geophysical, large-scale context on unstructured meshes.

Using a Stokeslet provides a mathematically consistent representation of pointwise interactions between icebergs and sea ice. The resulting drag term is incorporated into the sea-ice momentum equation, and the entire system is discretized using piecewise linear finite elements. We assume that each iceberg has a diameter smaller than the grid spacing. The approach avoids the need for mesh refinement and ensuring mesh-independent forcing. Our implementation is compatible with sea-ice components in ocean models that use unstructured meshes, such as FESOM and ICON \citep{Danilov2015, Mehlmann2021}.

To our knowledge, this is the first consistent numerical framework that captures the momentum exchange between individual, dynamically grounded subgrid icebergs and the sea-ice continuum. We further derive a stability estimate for the sea-ice momentum equation including the Stokeslet drag term.

The paper is structured as follows: Section \ref{sec:gov} introduces the governing model equations for the sea-ice and the iceberg model and derives an estimate for the viscous-plastic sea-ice model with subgrid iceberg--sea-ice drag term. Section~\ref{sec:disc} outlines the discretization of the coupled system. In Section \ref{sec:test_cases}, we numerically analyze the finite-element discretization and demonstrate that it is a suitable framework for fast-ice formation due to subgrid iceberg grounding. The paper ends with a conclusion in Section \ref{sec:con}.

\section{Governing equations}\label{sec:gov}
Let $\Omega \subseteq \mathbb{R}^2$ be  the two dimensional sea-ice domain of interest with a polygonal boundary and $I_t=[0,T]$ the time interval of which we consider the sea-ice evolution. By $x,y,t$ we denote the horizontal spatial coordinates and the time, respectively. The sea-ice motion is described by three variables: sea-ice concentration $a(x,y, t) \in [0,1]$ (the fraction of a grid cell covered with sea ice), mean sea-ice thickness $h(x,y, t) \in [0, \infty)$, and sea-ice velocity $\vt(x, t) \in \mathbb{R}^2$. In addition. we consider a set of disk-shaped icebergs represented by {{particles  $\{p\}$}}. Every iceberg is equipped with a radius $r_p$ and a height $h_p$, which can vary between the icebergs. 

\subsection{Sea-ice model}
The sea-ice velocity is determined from the momentum equation for all $(x,y,t) \in \Omega \times (0,T)$:
\begin{equation}\label{eq:mom}
    \rho h \partial_t \vt =f_{\text{c}}+f_{\text{sh}}+f_{\text{o}}+f_{\text{a}}+f_{\text{ib}}+f_{\text{r}},
\end{equation}
where $\rho$ is the sea-ice density. 
The right-hand side captures, respectively, the Coriolis force, the sea-surface slope, the atmospheric and oceanic drag, the iceberg drag, and internal stresses. 
The Coriolis force is given as
\begin{equation}
    f_{\text{c}}=-\rho h f k \times \vt,
\end{equation}
   with  Coriolis parameter $f$, the upward pointing unit vector $k$ and the gravity constant $g$. 
We approximate the force from the changing sea-surface height by
\begin{align}
    f_{\text{sh}} \approx \rho h f k \times \vt_{\text{o}},
\end{align} where $\vt_o$ is the ocean velocity \cite{Coon1980}.
The oceanic drag $f_{\text{o}}$ and atmospheric drag $f_{\text{a}}$ are modeled as
\begin{align*}
     f_{\text{o}}=C_{\text{o}} \rho_{\text{o}}\left|\vt-\vt_{\text{o}}\right|_{2} \left(\vt-\vt_{\text{o}}\right), \quad  
     f_{\text{a}}=C_{\text{a}} \rho_{\text{a}}\left |\vt_{\text{a}}\right|_{2} \vt_{\text{a}},
\end{align*}
where $C_o, C_a$ are the oceanic and atmospheric drag coefficients, $\rho_o,\, \rho_a$ are the ocean and atmosphere densities, $\vt_{\text{a}}$ is the velocity field of near surface atmospheric flow, respectively and  $|\cdot |_{2}$ denotes the Euclidean norm. 
To model the interaction between sea-ice and icebergs the usual momentum equation is modified by adding an iceberg drag:
\begin{equation}
  f_\text{ib}= a C_{\mathrm{i}} \rho_{\text{b}}\left|\vt_{b}-\vt\right|_{2} (\vt_b-\vt) \psi,
\end{equation}
where $\rho_\text{b}$ is the iceberg density, $C_\text{i}$ is the constant iceberg drag parameter. The indicator function  $\psi$ is given as
\begin{equation}
    \psi = \begin{cases} 1 \text{ if } (x-x_p)^2+(y-y_p)^2 < r_p^2, \\
0\text{ if } (x-x_p)^2+(y-y_p)^2 \geq r_p^2.  \end{cases}
    \label{eq:threshold}
\end{equation}
The internal stress are included to the momentum equation by 
\begin{equation}
    f_r=\div(\sigmat).
\end{equation}
The relationship between the internal sea-ice stresses, $\sigmat$, and the strain rates, 
\begin{equation}
   \epsilont=\frac{1}{2}(\nabla \vt+\nabla \vt^T)=\epsilont' + \frac{1}{2} \tr(\epsilont) I \in \mathbb{R}^{2 \times 2}, 
\end{equation}
are modeled by the viscous-plastic rheology \citep{Hibler1979}.
To make the notation more compact, we introduce the rheology in terms of the trace and deviatoric parts of the strain rate tensor \cite{Mehlmann2019}:
\begin{align*}
\sigmat := \frac{1}{2} \zeta \epsilont'(\vt) + \zeta \tr(\epsilont(\vt)) I - \frac{P}{2} I,
\end{align*}
where $I \in \mathbb{R}^{2 \times 2}$ is the identity matrix. The viscosity $\zeta$ and the ice strength $P$ are given by
\begin{align*}
\zeta &:= \frac{P}{2 \max(\Delta_P(\vt), \Delta_\text{min})}, & &\Delta_\text{min}: = 2 \times 10^{-9},\\ \Delta_P(\vt)&:=\sqrt{ \frac{1}{2} \epsilont' : \epsilont' + \tr(\epsilont)^2}, & &
 P:=h P^\star \exp(C(1 - a)),
\end{align*}
where the  $:$ in $\Delta_P(\vt)$ denotes the Frobenius inner product between tensors, $P^\star$ is the ice strength parameter and $C$  is the ice concentration parameter. To ensure a smooth transition between the viscous closure, $\Delta_\text{min}$, and plastic regimes, $\Delta_P(\vt)$, we follow \cite{Kreyscher2000} and replace $\max(\Delta_P(\vt), \Delta_\text{min})$ by 
\begin{equation}
\Delta(\dot{\boldsymbol{\epsilon}}) = \sqrt{\Delta_P(\dot{\boldsymbol{\epsilon}})^2 + \Delta_\text{min}^2}.
\label{eq:delta}
\end{equation}
All parameters of the problem are collected in Table~\ref{tab:parameters}. 
The advection of the sea-ice thickness, $h$, and the sea-ice concentration, $a$, are given  for all $(x,y,t) \in \Omega \times (0,T)$ by 
\begin{align}\label{eq:trans}
\partial_t a+\div(a \vt)=0, \quad 
\partial_t h+\div(h \vt)=0.
\end{align}
 The sea-ice system is completed by prescribing the boundary conditions as
\begin{align*}
    \vt=\vt^0, &\,\text{on} \quad \Omega \times \{t=0\},\\
    a=a^0, h=h^0, &\, \text{on} \quad \Omega \times \{t=0\},\\
    \vt=0 ,&\, \text{on} \quad \delta\Omega \times I_t,\\
    a=a^\text{in}, h=h^{in}, &\, \text{on} \quad \Gamma_\text{in} \times I_t,
\end{align*}
where  $\Gamma_\text{in}:=\{ (x,y) \in \delta\Omega | n\cdot \vt <0 \}$ denotes the part of the boundary with inflow. 
\begin{table}[h!]
\centering
\begin{tabular}{llr}
\hline
\textbf{Symbol} & \textbf{Description} & \textbf{Value} \\
\hline
$h_p$   & iceberg height                          & 200 m \\
$\rho$ & sea-ice density                      & 900 kg m$^{-3}$ \\
$\rho_o$ & ocean density                        & 1025 kg m$^{-3}$ \\
$\rho_a$ & atmosphere density                        & 1.3 kg m$^{-3}$ \\
$\rho_b$ & iceberg density                         & 900 kg m$^{-3}$ \\
$C_i$   & sea-ice coefficient of resistance    & 1 \\
$C_{a}$ & atmosphere drag coefficient          & $2.5 \times 10^{-4}$ \\
$C_{o}$ & ocean drag coefficient               & $5 \times 10^{-4}$ \\
$C_{va}$ & atmosphere coefficient of resistance&0.4\\
$C_{vo}$ & ocean  coefficient of resistance& 0.85\\
$f$ & Coriolis parameter & $1.46 \times 10^{-4}\ \mathrm{s^{-1}}$ \\
$P^*$ & ice strength parameter & $27.5 \times 10^3\ \mathrm{N\,m^{-2}}$ \\
$C$ & ice concentration parameter & $20$ \\
$\Delta_{min}$ &viscous closure&   $ 1 \times10 ^{-9} \mathrm{s^{-1}}$
\\
\hline
\end{tabular}
\caption{Physical parameters used in the sea-ice and iceberg model. \label{tab:parameters}}
\end{table}

\subsection{Iceberg model}
The position of an iceberg particle $\xt_p=(x_p,y_p)$ is related to its velocity $\vt_b$ by
\begin{equation}
   \vt_b= \partial_t \xt_p.
\end{equation}
Similarly to the sea-ice momentum equation, the momentum equation governing the iceberg dynamics (e.g.,~\cite{Lichey2001}) is given as 
\begin{equation}
    M \frac{d\vt_b}{dt} = \bm{F}_c + \bm{F}_{sh} +\bm{F}_o + \bm{F}_a +  \bm{F}_{si}, \label{eq:momentum}
\end{equation}
where $M=\rho_b h_p \pi r^2_p$ is the iceberg mass and $\vt_b$ is the iceberg velocity. 
The forces on the right-hand side of Eq.~\eqref{eq:momentum} are the Coriolis force, $\bm{F}_c = M f k \times \vt_b,$ the force created by the slope of the sea surface, $\bm{F_}{sh}=M f k\times \vt_o,$ and the oceanic and atmospheric drag, $\bm{F}_{o}$ and $\bm{F}_{a}$, respectively. 

The latter two drag terms, $\bm{F}_o$ and $\bm{F}_a$, are given as
\begin{align}
\bm{F}_{o} &= \pi r_p^2 \rho_o \left( C_{o} + \frac{h_p}{r_p} C_{vo} \right) |{\vt}_o - {\vt}_b|_2 ({\vt}_o - {\vt}_b), \\
\bm{F}_{a} &= \pi r_p^2 \rho_a \left( C_{a} + \frac{h_p}{r_p} C_{va} \right) |{\vt}_a|_2 {\vt}_a,
\end{align}
where $C_{vo}, C_{va}$ are the body drag coefficients. The vertical area of the iceberg in contact with air and water is given by $\pi r_p h_p$, whereas the horizontal contact  area with air and water is given by  $\pi r_p^2$.
The iceberg--sea-ice interaction is modeled similar to \cite{Lichey2001} as 
\begin{equation}
\bm{F}_{si} =a\rho C_i r_p h_p |\vt - \vt_b|_2 (\vt - \vt_b),
\end{equation}
where $C_i$ is the sea-ice coefficient of resistance and $\vt$ is the sea-ice velocity evaluated at the center of the disk-shaped iceberg $\xt_p$. The system is complete by prescribing an initial position and velocity of the iceberg $\xt(0)=\xt^ 0$  and $\vt_b=\vt^ 0_b$.

\subsection{A-priori stability estimate}
In this section, we derive an a-priori stability estimate for the norm of velocity fields that satisfy the sea-ice momentum equation with a subgrid-scale iceberg-sea-ice drag term. This estimate serves as a diagnostic tool in Section~\ref{sec:energyNum} to assess whether the discrete solution remains qualitatively consistent with the underlying continuous model.

The detailed derivation of the stability estimate is provided in Appendix A. All assumptions made in the analytical setting, such as constant sea-ice concentration and thickness, are consistently reflected in the numerical experiments discussed in Section~\ref{sec:energyNum}. In addition, we adopt a linearized form of the ocean drag term following \cite[Chapter 6.1.4]{Lepp2010} denoted by $\bar C_o$. Thus, the ocean drag reads as \begin{align}\label{eq:ocean_drag}
    \bar f_o= \rho_o \bar C_o (\vt-\vt_o).
\end{align}
The numerical analysis is conducted for the functional $\Phi(\vt)$, defined by
\begin{align}\label{eq:estimate}
     \Phi(\mathbf{v}) :=   \int_0^T  \sum^{N_p}_{i=1}aC_{i}\rho_{b}\|\vt\|_{B_{\xt_p}}^3\, dt,
\end{align}
where $N_p$ is the number of grounded icebergs and $\| \cdot \|_{B_{\xt_p}}$ denotes the $L^2$-norm over the domain covered by the iceberg with center $\xt_p$.
Based Theorem 1, provided in the Appendix A, the following estimate holds
\begin{align}
   \Phi(\mathbf{v}) &\leq \int^T_0 \frac{1}{2 \rho_o \bar C_o} \|\mathcal{R}\|^2 dt + \|\rho \mathbf{v}(0)\|^2, 
    \label{eq:evp_bound}
\end{align}  
where  $\mathcal{R} := f_{sh} + f_a+ \rho_o \bar C_{o} \mathbf{v}_o$ and $\| \cdot \|$ denotes the $L^2$-norm over $\Omega$. 
The estimate in~\eqref{eq:evp_bound} shows that the solution of the sea-ice momentum equation~\eqref{eq:mom}, when evaluated through a functional involving the iceberg--sea-ice drag term, is bounded by the problem data. Accordingly, a similar bound should be reflected by the finite-element discretization. This is studied in Section~\ref{sec:energyNum}.

\section{Numerical discretization}\label{sec:disc}
\subsection{Weak formulation of the momentum equation}
In the following, we derive the weak form of the momentum equation (\ref{eq:mom}) which is the basis for the  finite-element discretization outlined in Section \ref{sec:disc} and particular of interest for the approximation of the sub-grid iceberg--sea-ice drag. The Sobolev space \( H^1(\Omega)^2 \) comprises vector-valued functions defined on a two-dimensional domain \( \Omega \subseteq \mathbb{R}^2 \) with square-integrable first-order derivatives. The subspace \( H_0^1(\Omega) \) consists of functions in \( H^1(\Omega) \) whose trace vanishes on the boundary \( \Gamma \). Multiplying equation~\eqref{eq:mom} by a test function \(\phi \in V = (H_0^1(\Omega))^2 \) and integrating over the domain yields the following weak formulation:

\begin{equation}\label{eq:momWeak}
   \int_{\Omega} \rho h \partial_t \vt \cdot \phit \, dxdy=\int_\Omega (f_{\text{c}}+f_{\text{sh}}+f_{\text{o}}+f_{\text{a}}+f_{\text{ib}}+f_{\text{r}}) \cdot \phit \, dxdy.
\end{equation}

\subsection{Finite element discretization}
The sea-ice model is implemented using the open-source academic software library \textsc{Gascoigne}~\citep{Gascoigne}, which is based on quadrilateral meshes. The computational domain \( \Omega \) is discretized into a set of non-overlapping open quadrilaterals:
\begin{equation}
   \Omega_h = \bigcup_{i=1}^N K_i, \quad K_i \cap K_j = \emptyset \quad \text{for all } K_i \neq K_j \in \Omega_h.
\end{equation}
Velocity unknowns are defined at the mesh vertices, while tracer variables are located at the cell centers. The velocity field is approximated in space using piecewise linear finite elements.
The space of the piecewise linear elements on $\Omega_h$ is given as:
\begin{align}
    V_h=\{\phi \in C(\Omega), \forall K \in \Omega_h: \phi|_K \circ T_k \in Q\}, \quad Q:=\text{span}\{ x_iy_j,\, 0\leq i,j\leq 1\},
\end{align}
where by $T_K$ we describe the iso-parameteric reference map. To approximate the sea-ice velocity, we define the discrete velocity space as { \( V_h^{\mathbf{v}} = V \cap V_h \)}.  
Based on this space, the finite-element formulation for the velocity field reads
\begin{align*}
    \vt_h=\sum^N_{i=1}\vt^i\phit_i,
\end{align*}
where $\vt^i=(\vt^i_x,\vt^i_y)$ are the coefficients of the velocity vector. The discrete version of the momentum equation~\eqref{eq:mom} is given by 
\begin{align*}
\sum_K \sum^3_{i,j=1} \int_K  \rho h \partial_t \vt^i_h \cdot \phit_j\, dxdy=\int_K &(f^i_{c,h}+f^i_{sh,h}+f^i_{a,h}+f^i_{o,h})\cdot \phi^j\, dxdy\\+&\int_K f_{ib,h} \cdot \phit^j \, dxdy - \int_{K}\nabla \sigmat^{i}_h\cdot \nabla \phit^{j} \, dxdy.
\end{align*}
We describe the approximation of the subgrid-scale coupling between sea ice and grounded icebergs. All remaining terms are discretized using standard finite-element techniques; see, for example \citep{Mehlmann2019}. 

Let \( B_K(x, y) \) denote the circular subregion within a cell \( K \) that is covered by an iceberg. We assume that the iceberg diameter is smaller than the cell size. Under this assumption, the integral representing the drag induced by the iceberg simplifies to
\begin{equation}\label{eq:ib_disc}
 \begin{aligned}
   \int_K f_{ib,h} \cdot \phi \,dx dy &=  \int_K C_{\mathrm{i}} \rho_{\mathrm{b}}|\vt_b-\vt_h |_2(\vt_b-\vt_h) \cdot \phi \,dxdy\\
   &=  -C_{\mathrm{i}} \rho_{\mathrm{b}}\int_{B} |\vt_h|_2 \vt_h \cdot \phi\, dxdy \approx- C_{\mathrm{i}} \rho_{\mathrm{b}} \pi r^2_p|\vt_h(\xt_p)|_2 \vt_h(\xt_p)\cdot \phi(\xt_p).
\end{aligned}   
\end{equation}

\subsection{Discretization of the coupled iceberg--sea-ice system in time}
We discretize the coupled iceberg--sea-ice system in time using a splitting approach that separates the system into three main components: the sea-ice momentum equation, sea-ice advection, and iceberg advection.

First, the momentum equation (Eq.~\eqref{eq:mom}) is solved.
In the momentum equation the iceberg-induced drag term is treated for simplicity explicitly in time, while all other terms are handled implicitly. The implicit Euler method is commonly used to discretize the sea-ice momentum equation due to its stiffness, as explicit time discretizations demand extremely small time steps~\citep{ip_effect_of_rheology_1991}. The resulting nonlinear system is solved using a modified Newton method~\citep{MehlmannRichter2016newton}. After solving the momentum equation, the tracer variables in the transport equation (Eq.~\eqref{eq:trans}) are updated applying an upwind scheme. Once the sea-ice system has been computed, the iceberg model is advanced with an explicit Euler scheme. The complete time discretization of the coupled system is summarized in Algorithm~\ref{alg:dy}.

\begin{algorithm}
\caption{Time loop \label{alg:dy}}
Let \( I_t = [0, T] \) denote the time interval of interest, and let the initial conditions at time \( t_0 = 0 \) be given by the sea-ice velocity \( {\vt}(0)=\vt^0 \), sea-ice concentration \( a(t_0)=a^0 \), sea-ice thickness \( h(t_0)=h^0 \), as well as the positions \( \{ {\xt}_p(t_0)=\xt_p^0 \} \) and velocities \( {\vt}_b(t_0)=\vt_b^0 \) of the icebergs.
The time interval is divided into \( N \) equidistant steps: \\ 
\[
0 = t_0 < t_1 < \dots < t_N = T.
\]
The following time-stepping procedure is carried out for \( n = 1, 2, \dots, N \):
\begin{enumerate}
    \item Solve the sea-ice momentum equation (Eq.~\eqref{eq:mom}) using \( a(t_{n-1}) \) and \( h(t_{n-1}) \):
    \[
    {\vt}(t_{n-1}) \rightarrow {\vt}(t_n).
    \]
    
    \item Solve the advection equation (Eq.~\eqref{eq:trans}) using the velocity \( {\vt}(t_n) \):
    \[
    a(t_{n-1}) \rightarrow a(t_n), \quad h(t_{n-1}) \rightarrow h(t_n).
    \]
    
    \item Solve the iceberg momentum equation (Eq.~\eqref{eq:momentum}) using \( {\vt}_b(t_{n-1}) \), as well as \( {\vt}(t_n) \) and \( a(t_n) \):
    \[
   {\vt}_b(t_{n-1}) \rightarrow {\vt}_b(t_n).
    \]
    
    \item Update the iceberg positions using \( {\vt}_b(t_n) \) and \( \boldsymbol{x}_p(t_{n-1}) \):
    \[
    \boldsymbol{x}_p(t_{n-1}) \rightarrow \boldsymbol{x}_p(t_n).
    \]
\end{enumerate}
\end{algorithm}

\begin{figure}
    \centering
    \includegraphics[width=1.0\linewidth]{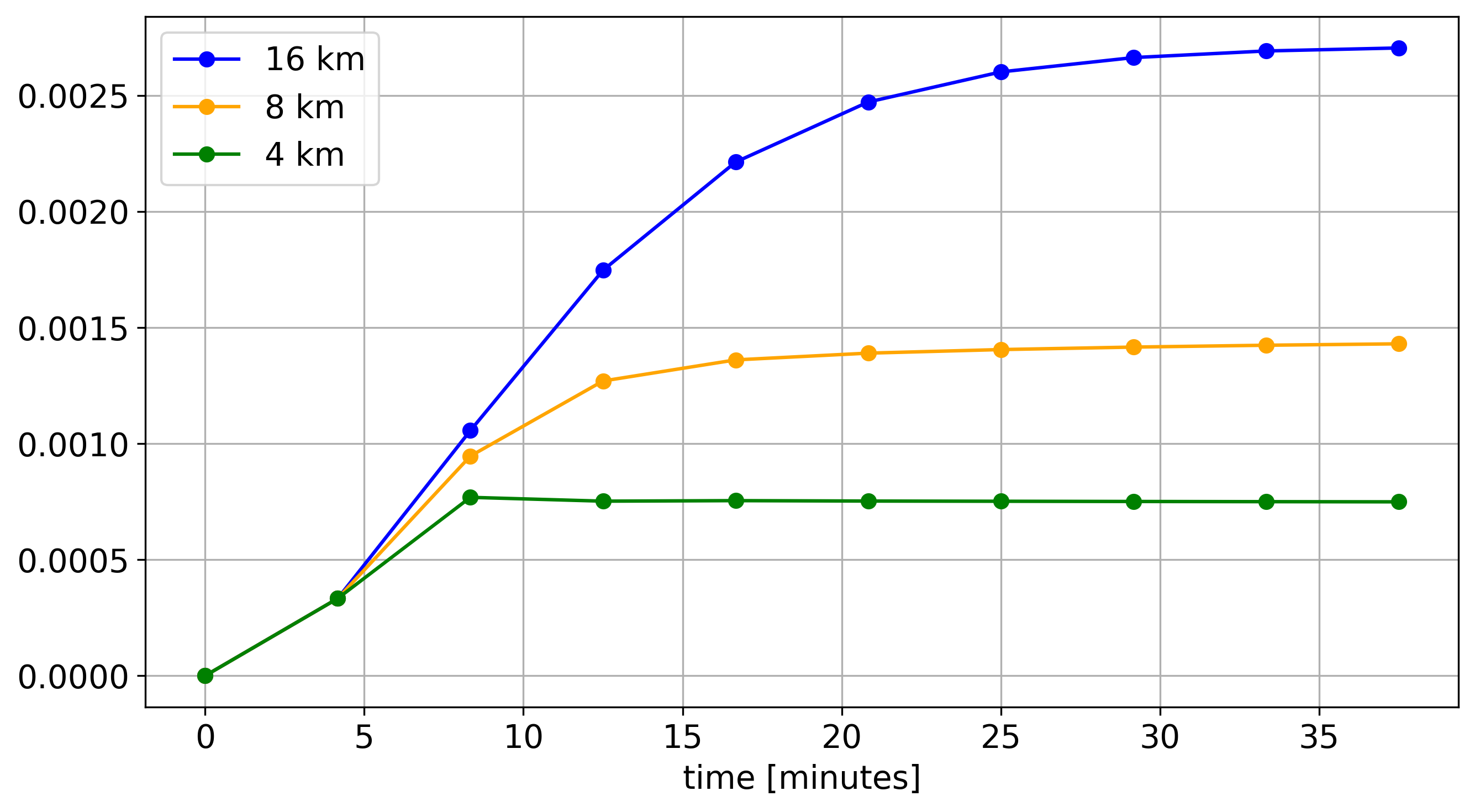}
    \caption{Estimate Eq.~\eqref{eq:estimate} evaluated over time for three different grid resolutions.}
    \label{fig:functional}
\end{figure}
\section{Numerical experiments}\label{sec:test_cases}
The viscous-plastic sea-ice model with a subgrid-scale iceberg--sea-ice drag term is analyzed in four different configurations. In the first setup (Section~\ref{sec:energyNum}) we numerically study the boundedness of the estimate introduced in~\eqref{eq:evp_bound}, to demonstrate that the approximated solution of the system is qualitatively consistent with the analytical solution.  In Section~\ref{sec:numA}, we investigate the influence of two subgrid-scale grounded icebergs on polynya formation and fast-ice development under grid refinement. Section~\ref{sec:SizeIceberg} presents an analysis of fast-ice formation in relation to grounded iceberg size.  Section~\ref{sec:numDyn} evaluates the full iceberg--sea-ice system, including traveling and dynamically grounding icebergs. 

All numerical experiments are performed on a \SI{512}{\kilo\meter} × \SI{512}{\kilo\meter} domain. The initial sea-ice concentration is set to 0.5, with an initial sea-ice thickness of \SI{1}{\meter}. Icebergs are modeled as \SI{200}{\meter}-high cylindrical bodies. Wind stress and the Coriolis force are neglected in this study. The ocean flows at a constant velocity of \SI{0.05}{\meter\per\second}, while both sea ice and icebergs are initially at rest.

\subsection{Stability estimate} \label{sec:energyNum}
To  evaluate the stability estimate, we switch-off the advection of the sea-ice thickness and sea-ice concentration. We place two grounded icebergs with a radius of $\SI{1}{\kilo\meter}$ at coordinates (\SI{159}{\kilo\meter}, \SI{159}{\kilo\meter}) and (\SI{159}{\kilo\meter}, \SI{157}{\kilo\meter}), respectively.
The estimate in~\eqref{eq:evp_bound} states that the functional $\Phi(\vt)$ derived for the solution to the continuous sea-ice equation with a subgrid iceberg--sea-ice drag term is bounded by the data. This property should also hold for the discrete approximations at different mesh levels. In Figure~\ref{fig:functional}, we evaluate the functional on a $\SI{16}{\kilo\meter}$, $\SI{8}{\kilo\meter}$, and $\SI{4}{\kilo\meter}$ mesh. The functional remains bounded and converges toward a stationary solution, as expected, since both the forcing and the sea-ice tracers are held constant over time. We observe that the difference between the evaluated functionals decreases with increasing spatial resolution, indicating convergence in space.  



\subsection{Fast-ice formation on different gird levels}\label{sec:numA} 
\begin{figure}
    \centering
    \begin{tabular}{ c c c c}
        \includegraphics[scale=0.5]{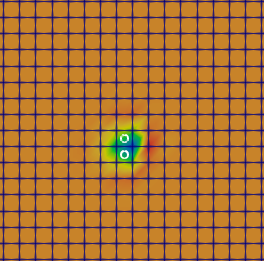} & 
        \includegraphics[scale=0.5]{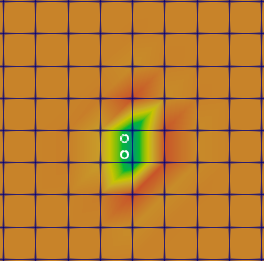} & 
        \includegraphics[scale=0.5]{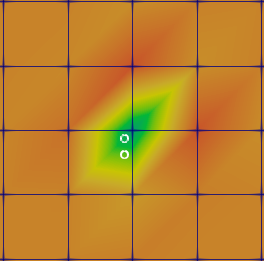} & 
        \includegraphics[scale=0.4]{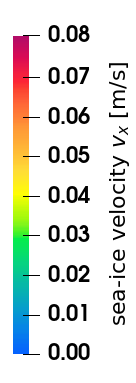} \\
        \includegraphics[scale=0.4]{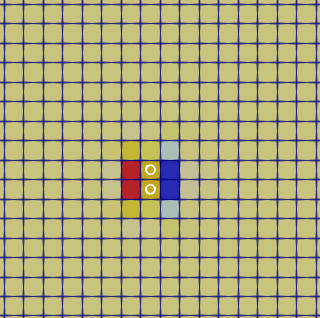} & 
        \includegraphics[scale=0.4]{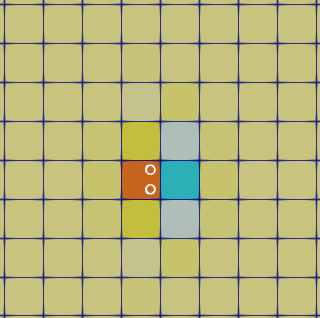} & 
        \includegraphics[scale=0.4]{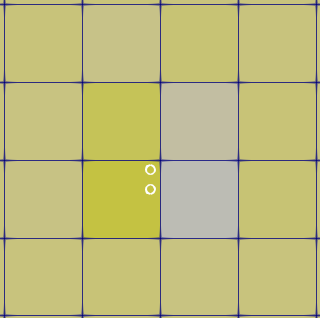} & 
        \includegraphics[scale=0.39]{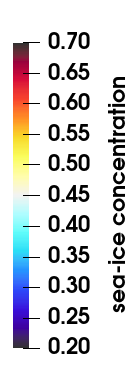} \\
        \SI{4}{\kilo\meter} & \SI{8}{\kilo\meter} & \SI{16}{\kilo\meter} & \\
    \end{tabular}
    \caption{Closeup of the sea-ice velocity ($\vt_x$) in the upper row and sea-ice concentration in the lower row for varying grid resolutions. All closeups show the same area ($\SI{128}{\kilo\meter} \leq (x, y) \leq \SI{192}{\kilo\meter}$). The grounded iceberg particles are marked by white circles.}
    \label{fig:tc1_vel}
\end{figure}

 We place the center of two  grounded  iceberg  with a radius of \SI{1}{\kilo\meter} at  coordinates {$(\SI{158}{\kilo\meter} ,\SI{158}{\kilo\meter})$ and $(\SI{158}{\kilo\meter}, \SI{154}{\kilo\meter} )$}, respectively  and analyze their feedback on the sea-ice model. In the upper two panels of Figure~\ref{fig:tc1_vel}, we show the first component of the sea-ice velocity field. The influence of the subgrid-scale grounded iceberg is clearly visible: the velocity drops near the nodes (i.e., the vertex of a quadrilateral cell) that lies closest to the position of grounded iceberg. In the case of the $\SI{16}{\kilo\meter}$ mesh (left panel in Figure \ref{fig:tc1_vel}), the icebergs address the upper right node, as both icebergs ground in the upper right region of this cell.

 On the $\SI{8}{\kilo\meter}$ mesh (middle panel in Figure \ref{fig:tc1_vel}), the influence of the grounded iceberg is visible on the upper right and lower right node. Due to the refinement of the mesh from $16\mathrm{km}$ to $8\mathrm{km}$ the icebergs are now located near the midpoint between these two nodes. Increasing the resolution to $\SI{4}{\kilo\meter}$ (right panel in Figure \ref{fig:tc1_vel}) assigns the icebergs to two distinct cells. In this case, the upper iceberg affects the lower right node of its cell, while the lower iceberg influences the upper right node of its respective cell.  This behavior is consistent with the finite-element interpolation used in \eqref{eq:ib_disc}. 
 \begin{figure}
    \centering
    \includegraphics[width=1.0\linewidth]{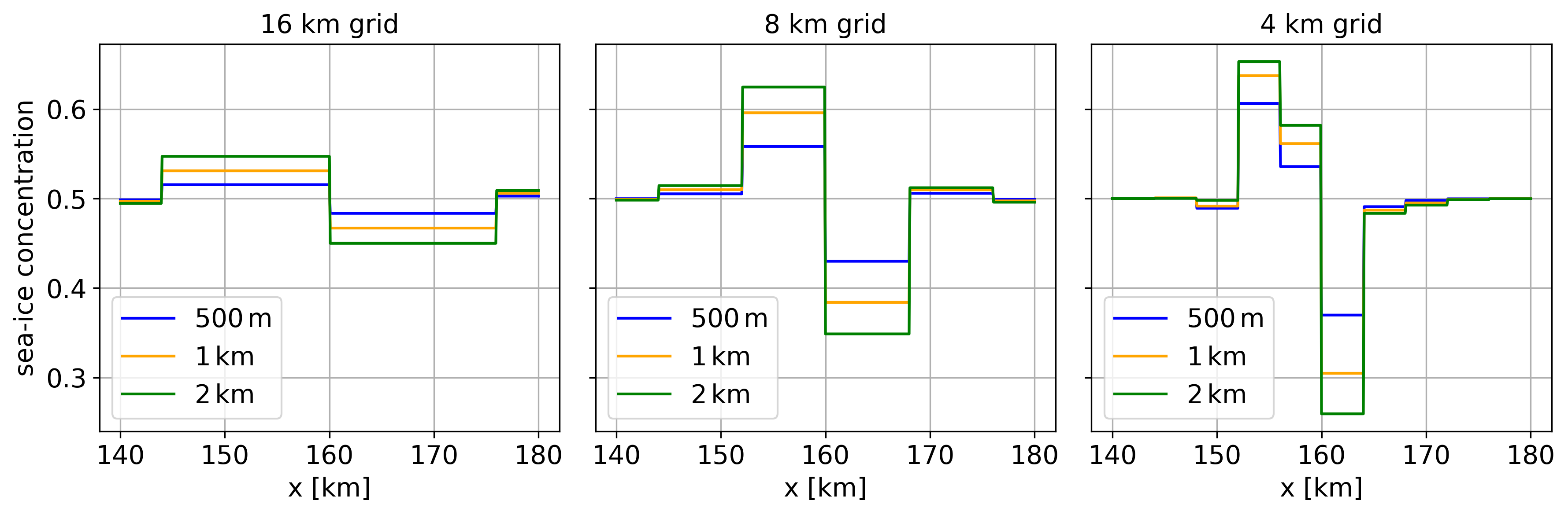}
    \caption{Cross-section of sea-ice concentration along $y = \SI{158}{\kilo\meter}$, intersecting the center of the upper iceberg $(\SI{158}{\kilo\meter} ,\SI{158}{\kilo\meter})$.  Panels show the ice concentration for grid sizes of \SI{16}{\kilo\meter}, \SI{8}{\kilo\meter}, and \SI{4}{\kilo\meter}, each with iceberg radii of \SI{500}{\meter}, \SI{1}{\kilo\meter}, and \SI{2}{\kilo\meter}.
   }
    \label{fig:plotoverline}
\end{figure}

\begin{figure}
    \centering
    \includegraphics[width=1.0\linewidth]{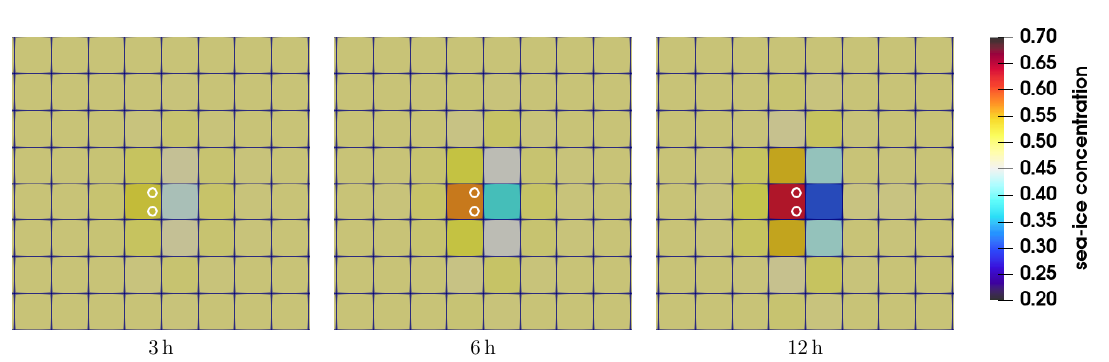}
    \caption{Time evolution of polynya formation. We show the sea-ice concentration in the region $\SI{128}{\kilo\meter} \leq (x, y) \leq \SI{192}{\kilo\meter}$ on an \SI{8}{\kilo\meter} grid. Grounded icebergs are indicated by white circles.}
    \label{fig:DynPol}
\end{figure}
In Figure~\ref{fig:tc1_vel} in the lower row, we show  the corresponding sea-ice concentration under grid refinement. On all considered mesh levels, the grounded icebergs lead to an increase in sea-ice concentration and thickness (not shown) in the grid cells to the left of the grounding position and a decrease in those to the right. The decrease can be associated with a polynya opening, while the increase can be attributed to the formation of fast ice. 
The intensity of the formed polynya as well as the number of affected grid cells depends on the mesh resolution, which in turn determines the influence of grounded icebergs on the neighboring velocity nodes, see upper panels in Figure~\ref{fig:tc1_vel}. By comparing the cross-sectional plot along the grounding position of the upper iceberg in Figure~\ref{fig:plotoverline}, we observe that a higher spatial resolution leads to narrower but more pronounced polynyas.

To further investigate the influence of the grounded iceberg on the polynya development and the fast-ice formation, we analyze the temporal dynamics of the setup. In Figure~\ref{fig:DynPol}, we visualize the evolution of the polynya over time on the  $\SI{4}{\kilo\meter}$ mesh. The polynya opens rapidly and becomes increasingly pronounced as the simulation progresses. The fast ice which develops up-stream to the polynya stays stable over the simulation period of 10 days (not shown). 
\begin{figure}
    \centering
    \begin{tabular}{ c c }
        \includegraphics[scale=0.5]{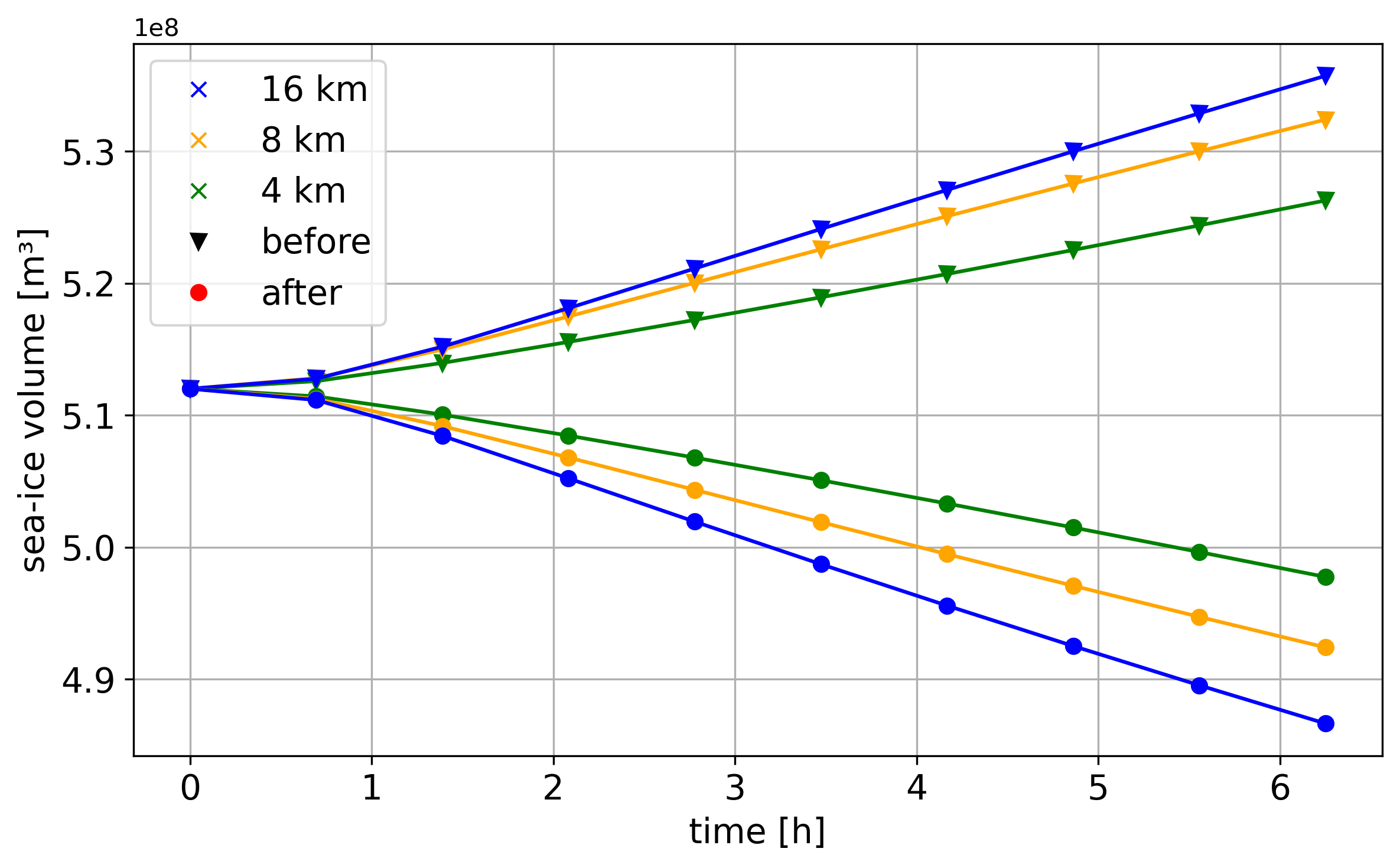} & 
        \includegraphics[scale=0.5]{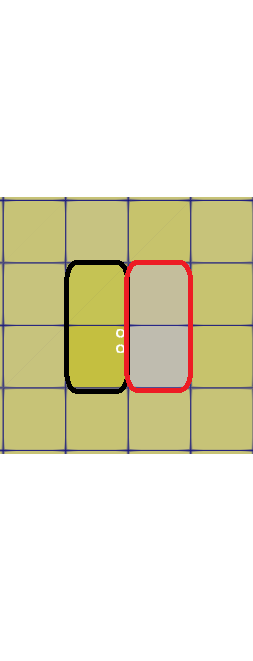} \\
    \end{tabular}
    \caption{
    Sea-ice volume over time in iceberg-affected cells for three grid resolutions. Triangles and circles indicate integration over the black ("before") and red ("after") marked area, respectively. Right: closeup of the area $\SI{128}{\kilo\meter} \leq (x, y) \leq \SI{192}{\kilo\meter}$ at \SI{16}{\kilo\meter} resolution.
    }
    \label{fig:InCon}
\end{figure}
To analyze the dynamics across different mesh resolutions, we integrate the sea-ice concentration over a rectangular region located before and after the grounding position,  see Figure~\ref{fig:InCon}. We observe  a converging behavior of the sea-ice concentration in the area before and after the grounding position as the  difference between the integral value decreases with increasing resolution.

\subsection{Fast-ice formation due to different iceberg sizes}\label{sec:SizeIceberg}
\begin{figure}
    \centering
    \includegraphics[width=1.0\linewidth]{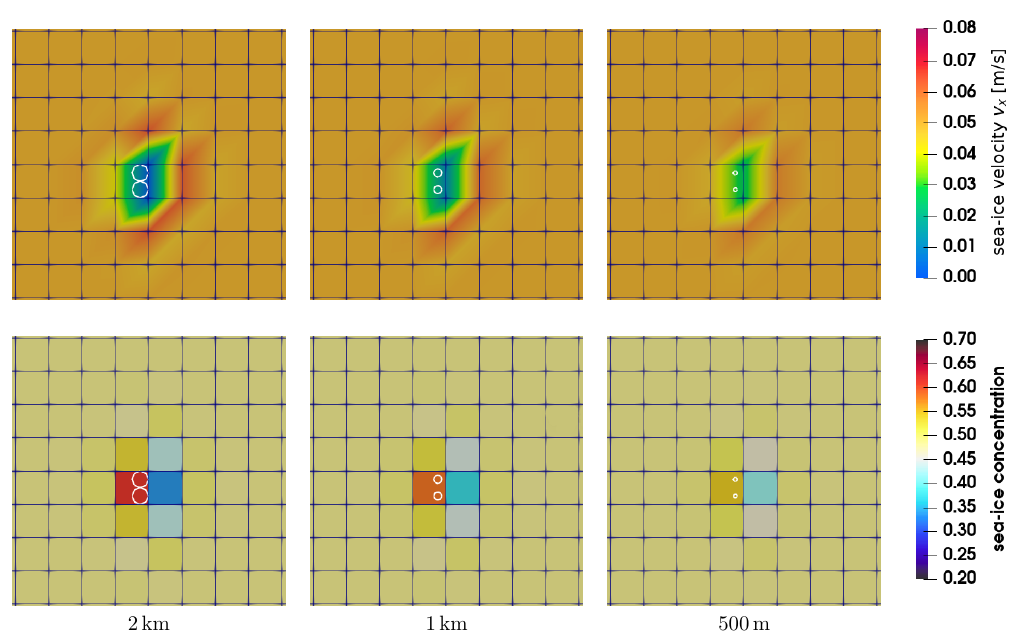}
     \caption{Closeup of the sea-ice velocity ($\vt_x$) in the upper row and sea-ice concentration in the lower row for different iceberg radii. The grounded iceberg particles are marked by white circles. All closeups show the same area ($\SI{128}{\kilo\meter} \leq (x, y) \leq \SI{192}{\kilo\meter}$) on a \SI{8}{\kilo\meter} grid.}
     \label{fig:bergRadius}
 \end{figure}

As in the previous example, we position two icebergs at the coordinates (\SI{158}{\kilo\meter}, \SI{158}{\kilo\meter}) and $(\SI{158}{\kilo\meter}, \SI{154}{\kilo\meter})$, respectively. In this setup, we fix the grid resolution to $\SI{8}{\kilo\meter}$  and systematically reduce the iceberg radius from $\SI{2}{\kilo\meter}$ and $\SI{1}{\kilo\meter}$ down to $\SI{0.5}{\kilo\meter}$. The x-component of the velocity filed is presented in the upper panel of  Figure~\ref{fig:bergRadius}. We observe that influence of the iceberg--sea-ice drag decreases with reduced iceberg radius. This is expected, as the drag term is weighted according to the size of the respective iceberg, see Equation~\eqref{eq:ib_disc}.
The lower panel of Figure~\ref{fig:bergRadius} shows the resulting sea-ice concentration fields. A reduction in sea-ice concentration appears in the downstream region behind each iceberg, while sea ice accumulates upstream (to the left of the iceberg). As the iceberg radius decreases, this influence becomes progressively weaker. The cross-sectional plot presented in Figure~\ref{fig:InCon} indicates that this observation holds true for the 
$\SI{16}{\kilo\meter}$, 
$\SI{8}{\kilo\meter}$ and 
$\SI{4}{\kilo\meter}$ mesh.

\subsection{Dynamical iceberg grounding and traveling icebergs} \label{sec:numDyn}
\begin{figure}
    \centering
    \makebox[\columnwidth]{
    \begin{tabular}{ c c c c}
        \includegraphics[scale=0.3]{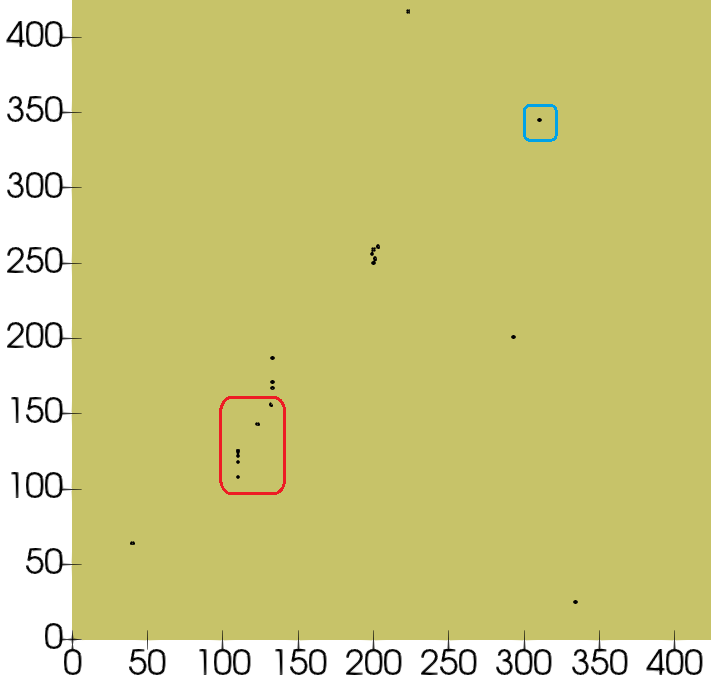} & 
        \includegraphics[scale=0.3]{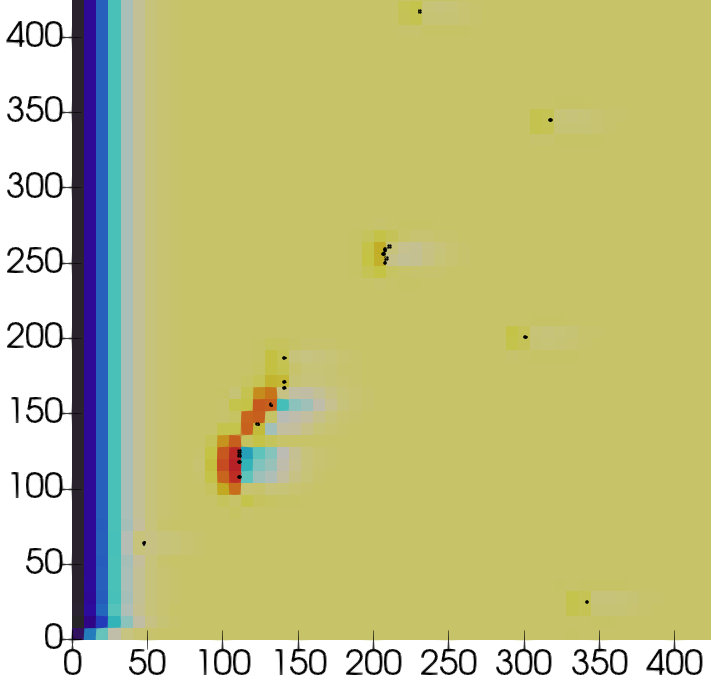} &
        \includegraphics[scale=0.4]{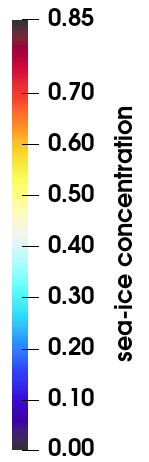}\\
       & 0 days & 3 days & \\
    \end{tabular}
    }
    \caption{Sea-ice concentration of a scenario with multiple icebergs. The red box on the left marks the icebergs which will ground during the simulation. The blue box marks an individual iceberg, which will not ground.}
    \label{fig:complex_tc}
\end{figure}

We consider a scenario involving 19 dynamic icebergs, each with a radius of \SI{1}{\kilo\meter} on a mesh with  \SI{8}{\kilo\meter} resolution. The icebergs are placed at various locations, as listed in Table~\ref{tab:split_not_grounded}. The initial setup is visualized in the right panel of Figure~\ref{fig:complex_tc}. In this test case, we simulate a dynamic grounding event and analyze the effect of icebergs which are traveling through sea-ice on a subgrid scale. 

Grounding events typically occur in shallow waters and have a significant impact on sea-ice dynamics. When icebergs come into contact with the seafloor, they become immobilized. We assume the seafloor to be shallow within a rectangular area defined by the lower left corner at (\SI{111}{\kilo\meter}, \SI{100}{\kilo\meter}) and the upper right corner at (\SI{200}{\kilo\meter}, \SI{165}{\kilo\meter}). Icebergs that ground during the simulation are highlighted with a red box in Figure~\ref{fig:complex_tc}.

The right panel of Figure~\ref{fig:complex_tc} shows the resulting sea-ice concentration after three days of simulation. The icebergs entering the shallow area are grounded (we set $\vt_b=0$) and contribute to the formation of a polynya downstream of the grounding position, as well as to the development of fast ice upstream. Even icebergs that do not ground influence the surrounding sea-ice concentration.

To illustrate the effect of moving icebergs in detail, we analyze the influence of a single iceberg on the surrounding sea-ice cover. The initial position of the iceberg is marked by a blue box in Figure~\ref{fig:complex_tc}, with its center located at coordinates (\SI{310}{\kilo\meter}, \SI{345}{\kilo\meter}).

\begin{figure}
    \centering
        \includegraphics[width=0.8\linewidth]{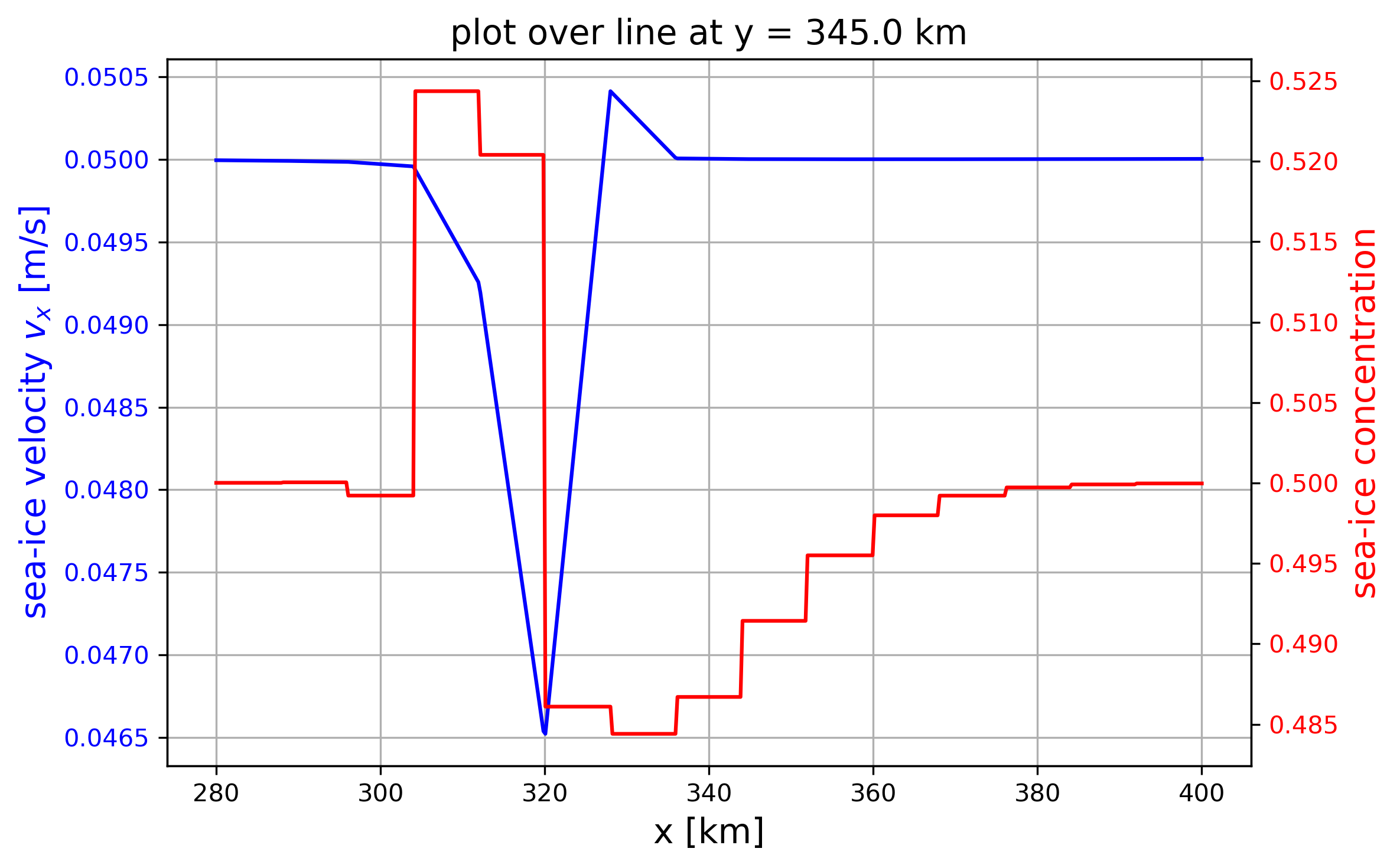}
    \caption{Cross-section of the sea-ice concentration (red) and velocity $\vt_x$ (blue) along $y = \SI{345}{\kilo\meter}$, intersecting the ungrounded iceberg at (\SI{310}{\kilo\meter}, \SI{345}{\kilo\meter}).}
    \label{fig:pol_complex}
\end{figure}

The cross-sectional plot in Figure~\ref{fig:pol_complex} shows the sea-ice concentration after 3 simulated days. We observe a slight accumulation of sea ice to the left of the iceberg, along with a local decrease in concentration to the right (red line). This pattern results from the iceberg moving slightly more slowly than the surrounding sea ice. As a consequence, the iceberg--sea-ice drag term~\eqref{eq:ib_disc} locally reduces the sea-ice velocity near the iceberg's position, as indicated by the blue line in Figure~\ref{fig:pol_complex}. This local deceleration of sea ice leads to an upstream accumulation (left of the grounding position) and a downstream reduction in sea-ice concentration.

\begin{table}
    \centering
    \footnotesize
    \begin{tabular}{p{4.0cm}|p{4.0cm}|p{4.0cm}}
    \textbf{Grounded} & \textbf{Not grounded} & \textbf{Not grounded } \\
    \hline
    \shortstack[l]{
        (\SI{110}{\kilo\meter}, \SI{108}{\kilo\meter})\\
        (\SI{110}{\kilo\meter}, \SI{118}{\kilo\meter})\\
        (\SI{110}{\kilo\meter}, \SI{122}{\kilo\meter})\\
        (\SI{110}{\kilo\meter}, \SI{125}{\kilo\meter})\\
        (\SI{123}{\kilo\meter}, \SI{143}{\kilo\meter})\\
        (\SI{132}{\kilo\meter}, \SI{156}{\kilo\meter})
    } &
    \shortstack[l]{
        (\SI{40}{\kilo\meter}, \SI{64}{\kilo\meter})\\
        (\SI{133}{\kilo\meter}, \SI{167}{\kilo\meter})\\
        (\SI{133}{\kilo\meter}, \SI{171}{\kilo\meter})\\
        (\SI{133}{\kilo\meter}, \SI{187}{\kilo\meter})\\
        (\SI{199}{\kilo\meter}, \SI{256}{\kilo\meter})\\
        (\SI{200}{\kilo\meter}, \SI{250}{\kilo\meter})\\
        (\SI{200}{\kilo\meter}, \SI{259}{\kilo\meter})
    } &
    \shortstack[l]{
        (\SI{201}{\kilo\meter}, \SI{253}{\kilo\meter})\\
        (\SI{203}{\kilo\meter}, \SI{261}{\kilo\meter})\\
        (\SI{223}{\kilo\meter}, \SI{417}{\kilo\meter})\\
        (\SI{293}{\kilo\meter}, \SI{201}{\kilo\meter})\\
        (\SI{310}{\kilo\meter}, \SI{345}{\kilo\meter})\\
        (\SI{334}{\kilo\meter}, \SI{25}{\kilo\meter})
    } \\
    \end{tabular}
    \caption{Initial position of the iceberg centers. The label grounded and not grounded indicates if the iceberg become grounded during the simulation. }
    \label{tab:split_not_grounded}
\end{table}

\section{Conclusion}\label{sec:con}
In this paper, we introduced a novel subgrid-scale drag formulation to represent the interaction between icebergs  and sea ice on coarse horizontal meshes. The drag is modeled using a Stokeslet, which provides a localized, mathematically consistent force representation at subgrid scale. The coupled system is discretized  with linear finite elements.
To establish theoretical consistency, we derived a stability estimate that bounds a functional of the sea-ice system by the problem data. Numerical evaluation of the functional $\Phi(\mathbf{v})$ confirms that the modified system remains bounded in the presence of subgrid-scale forcing.

We demonstrated that the inclusion of the Stokeslet-based iceberg--sea-ice drag term enables the simulation of fast-ice formation  and polynya opening caused by icebergs grounded at subgrid scale.  The amount of simulated fast ice and the intensity of polynya openings depend on the area covered by grounded icebergs and the chosen grid resolution. We further analyzed the impact of moving icebergs on the sea-ice cover and found that, when icebergs and sea-ice move in the same direction, but the iceberg velocities are lower than those of the surrounding sea ice, the drag term causes a slight accumulation of sea ice upstream (behind the iceberg) and a reduction downstream.

Future work will extend the method to realistic Antarctic configurations and evaluate its performance against observational data.

\section*{A \quad Appendix}

\begin{theorem}[A priori stability estimate]  Let  \( a \) and \( h \) be constant in time and space and \( a, h > 0 \). Assume \( \mathbf{v} \in L^2(I_t, \Vt_h^\vt) \) is a solution of the sea-ice momentum equation~\eqref{eq:momWeak} with $N_p$ grounded icebergs ($\vt_p=0$). Then, \( \mathbf{v} \) satisfies the following estimate
\end{theorem}
\begin{align}
\|\rho h \mathbf{v}(T)\|^2 &+ \int_0^T \{ \rho_o \frac{\bar C_{o}}{2} \|\mathbf{v}\|^2+ 
\|\sqrt{2}^{-1} \sqrt{\zeta} \epsilont \|^2 + \frac{3}{4} \|  \sqrt{\zeta} \operatorname{tr}(\epsilont) \|^2 \,
dt\\
&+\int_{0}^T\sum^{N_p}_{i=1}aC_{i}\rho_{b}\|\vt\|_{B_{\xt_p}}^3\,dt
\leq \int_0^T \frac{1}{2 \rho_o \bar C_o} \|\mathcal{R}\|^2 dt + \|\rho h \mathbf{v}(0)\|^2,
\end{align}
where \( \mathcal{R} := f_{sh} + f_a+ \rho_o \bar C_{o} \mathbf{v}_o\), $N_p$ is the number of grounded icebergs and $\| \cdot \|_{B_{\xt_p}}$ denotes the $L^2$-norm over the domain covered by the iceberg with center $\xt_p$.

\begin{proof}
Note that a $H^1$-estimate for the sea-ice model without the iceberg--sea-ice drag was carried out in~\citep{Mehlmann2021}.  We begin with the weak formulation of the momentum equation~\eqref{eq:momWeak} and select the test function \( \phi = {\vt} \).
\begin{equation}\label{eq:weak_energy}
   \int_{\Omega} \rho h \partial_t \vt \cdot \vt \, dxdy -\int_\Omega (f_{\text{c}}+f_{\text{o}}+f_{\text{ib}}+f_{\text{r}})\cdot \vt \, dxdy=\int_\Omega (f_{\text{sh}}+f_{\text{a}})\cdot \vt \, dxdy.
\end{equation}
We examine each term on the left-hand side of~\eqref{eq:weak_energy} individually: By applying the chain rule to the time-dependent integral we get
\begin{equation}\label{eq:ft}
\int_{\Omega} (\rho h \partial_t \mathbf{v} \cdot \mathbf{v}) \, dxdy = \frac{1}{2} \int_{\Omega} \rho h \partial_t |\mathbf{v}|_2^2 \, dxdy 
= \frac{1}{2} \partial_t \|\rho h \mathbf{v}\|^2.
\end{equation}
The integral of the Coriolis force vanishes because \( f_c \cdot {\vt} = (\rho h f \, \mathbf{k} \times {\vt}) \cdot {\vt} \) is antisymmetric. Assuming linear drag~\eqref{eq:ocean_drag}, the ocean drag simplifies to
\begin{align}\label{eq:fo}
\int_{\Omega}
\bar f_o \cdot \mathbf{v} \, dxdy 
= \int_\Omega\rho_o \bar C_{o} \mathbf{v}_o \cdot  \mathbf{v}- \rho_o\bar C_{o} \mathbf{v} \cdot \vt \, dxdy =:\int_\Omega \bar f_o^w-\bar f_o^\vt \, dx dy.
\end{align}
Following \citep{Mehlmann2021} the integral over the rheology can be written as 
\begin{align}\label{eq:fr}
  \int_{\Omega}-f_r\cdot \vt, dxdy=    \int_{\Omega}-\operatorname{div}(\boldsymbol{\sigma})\cdot \mathbf{v} \,dxdy  
     &= \|\sqrt{2}^{-1} \sqrt{\zeta} \epsilont \|^2 + \frac{3}{4} \|  \sqrt{\zeta} \operatorname{tr}(\epsilont) \|^2. 
\end{align}
Under the assumption icebergs  are  grounded ($\vt_b=0$), the integral over  $f_{ib}$ reduce to
\begin{equation}\label{eq:fib}
   -\int_{\Omega} f_{ib}\cdot \vt\, dxdy= \int_\Omega aC_{i}\rho_{b}|\vt|_2\vt \cdot \vt\, dxdy= \sum^{N_p}_{i=1} aC_{i}\rho^i_{b}\|\vt^i\|_{B_{\xt_p}}^3.
\end{equation}
To enhance clarity, we group together all terms that do not depend on $\vt$:
\[
\mathcal{R} := f_{sh}+f_a+ \bar f_o^w.
\]
Based on \eqref{eq:weak_energy} - \eqref{eq:fib}  we get 
\begin{equation}
\begin{aligned}\label{eq:lefthandsie}
\frac{1}{2} \partial_t \|\rho h \mathbf{v}\|^2+& \rho_o\bar C_{o}\|\vt\|^2
+ \sum^{N_p}_{i=1}\pi r^2_paC_{i}\rho_{b}{\|\vt\|^3 }\\&+ \|\sqrt{2}^{-1} \sqrt{\zeta} \dot{\boldsymbol{\varepsilon}} \|^2 + \frac{3}{4} \|  \sqrt{\zeta} \operatorname{tr}(\dot{\boldsymbol{\varepsilon}}) \|^2 
 \leq \int_\Omega \mathcal{R} \cdot \vt\, dxdy.
\end{aligned}
\end{equation}
By applying Youngs and Cauchy-Schwarz inequality  it holds that
\begin{align}\label{eq:rhs}
(\mathcal{R}, \mathbf{v}) 
&\leq \frac{1}{2\epsilon} \|\mathcal{R}\|^2 + \frac{\epsilon}{2} \|\mathbf{v}\|^2.
\end{align}
Choosing $\epsilon = \rho_o\bar C_o$ and integrating over $I_t$ yields the final estimate
 \begin{align*}
\|\rho h \mathbf{v}(T)\|^2& - \|\rho h \mathbf{v}(0)\|^2 
+ \int_0^T \{  \rho_o \frac{\bar C_{o}}{2}\|\mathbf{v}\|^2 
+\|\sqrt{2}^{-1} \sqrt{\zeta} \dot{\boldsymbol{\varepsilon}} \|^2\\
& + \frac{3}{4} \|  \sqrt{\zeta} \operatorname{tr}(\dot{\boldsymbol{\varepsilon}}) \|^2
+\sum^{N_p}_{i=1}\pi r^2_paC_{i}\rho_{b}\|\vt\|_{B_{\xt_p}}^3  \} dt
\leq \int_0^T \frac{1}{2\rho\bar C_o}\|\mathcal{R}\|^2 \, dt.
\end{align*}
\end{proof}

\bibliographystyle{plain}

\end{document}